%% file: tree17.tex
\documentclass{article}

\input{epsf.tex}
\usepackage{amsfonts,amssymb,verbatim,amsmath,amsthm,latexsym,textcomp,amscd}
\usepackage{latexsym,amsfonts,amssymb,epsfig,verbatim}
\usepackage{amsmath,amsthm,amssymb,latexsym,graphics,textcomp}
\usepackage{graphicx}
\usepackage{color}
\usepackage{url}

\input xy
\xyoption{all}

\theoremstyle{plain}
\newtheorem{theorem}{Theorem}[section]

\newtheorem{lemma}[theorem]{Lemma}
\newtheorem{corollary}[theorem]{Corollary}

\newtheorem{proposition}[theorem]{Proposition}

\newtheorem*{corollary*}{Corollary}

\newtheorem*{theorem_mt1*}{Theorem \ref{T:main tree 1}}
\newtheorem*{theorem_cc1*}{Theorem \ref{T:convex ccpt 1}}

\def\Mod{\mbox{\rm{Mod}}}

\def\Diff{\mbox{\rm{Diff}}}
\def\Stab{\mbox{\rm{Stab}}}

\def\Aut{\mbox{\rm{Aut}}}
\def\Out{\mbox{\rm{Out}}}
\def\Inn{\mbox{\rm{Inn}}}

\def\diam{\mbox{\rm{diam}}}

\def\intr{\mbox{\rm{int}}}
\def\d{\mbox{\rm{d}}}

\def\C{{\mathcal C}}
\def\D{{\mathcal D}}

\def\H{{\mathcal H}}

\def\hull{\mbox{\rm{Hull}}}

\def\co{\colon\thinspace}

\newcommand{\Id}{\operatorname{Id}}

\title{\textbf{Trees and mapping class groups}}
\author{Richard P. Kent IV\thanks{This research was partially conducted during the period the first author was employed by the Clay Mathematics Institute as a Liftoff Fellow.  This author was also supported by a Donald D. Harrington Dissertation Fellowship and an NSF postdoctoral fellowship.}, Christopher J. Leininger\thanks{partially supported by NSF grant DMS-0603881}, Saul Schleimer\thanks{partially supported by NSF grant DMS-0508971}}

\begin{document}
\maketitle

\footnotetext{This work is in the public domain.}

\begin{abstract}
There is a forgetful map from the mapping class group of a punctured
surface to that of the surface with one fewer puncture.
We prove that finitely generated purely pseudo-Anosov subgroups of the kernel of 
this map are convex cocompact in the sense of B. Farb and L. Mosher.
In particular, we obtain an affirmative answer to
their question of local convex cocompactness of K. Whittlesey's group.

In the course of the proof, we obtain a new proof of a theorem of I. Kra.
We also relate the action of this kernel on the curve complex to a family of 
actions on trees.
This quickly yields a new proof of a theorem of J. Harer.

\end{abstract}

\section{Introduction} \label{S:intro}

Let $S = S_{g,m}$ be an orientable surface of genus $g$ with $m$ punctures and assume throughout that the complexity
$\xi(S) = 3g-3+m$ is at least $1$.
Let $(S,z)$ denote the surface $S$ equipped with a marked point $z$.
There
is a homomorphism from the mapping class group of $(S,z)$
to that of $S$ which
fits into J. Birman's exact sequence \cite{birmansequence,birmanbook}
\[
1 \to \pi_1(S,z) \to \Mod(S,z) \to \Mod(S) \to 1,
\]
see \S \ref{S:background}.
We use this sequence to view $\pi_1(S,z)$ as a subgroup
of $\Mod(S,z)$.

Our main theorem answers Question 6 of \cite{kentleininger}.
\begin{theorem_cc1*}
If $G < \pi_1(S,z) < \Mod(S,z)$ is finitely generated and purely pseudo-Anosov, then $G$ is convex cocompact.
\end{theorem_cc1*}

Convex cocompactness for subgroups of the mapping class group was defined by B. Farb and L. Mosher in \cite{FMcc} via
the action on Teichm\"uller space by way of analogy with Kleinian groups and their action on hyperbolic space: a subgroup is convex cocompact if it has a quasiconvex orbit in the Teichm\"uller space.
Their work exhibits an intimate connection between convex cocompactness of a subgroup of the mapping class group and the geometry of the associated
surface group extension.
This concept was studied further by the first two authors in \cite{kentleininger}, extending the analogy with Kleinian groups; and by U. Hamenst\"adt in \cite{hamenstadt}, where the connection to surface group
extensions was strengthened. We note that Theorem \ref{T:convex ccpt 1} provides the first nontrivial examples of
convex cocompact groups that do not arise from a combination or ping--pong argument.

In terms of the analogy with Kleinian groups, Theorem \ref{T:convex ccpt 1} should be compared with a theorem of G. P. Scott
and G. A. Swarup \cite{scottswarup}:  finitely generated subgroups of infinite index in fiber subgroups of fibered hyperbolic
$3$--manifold groups are geometrically finite.
Indeed, the exact sequence of such a fibration
\[
1 \to \pi_1(S,z) \to \pi_1(M,z) \to \mathbb Z \to 1
\]
injects into Birman's sequence, and the subgroups covered by Theorem \ref{T:convex ccpt 1} are natural analogues
of those considered by Scott and Swarup.

Convex cocompact groups are necessarily finitely generated and virtually purely pseudo-Anosov \cite{FMcc}.  An important
question is whether the converse holds---this is Question 1.5 of \cite{FMcc} (see also Problem 3.4 of
\cite{Mosherproblems}), asked by Farb and Mosher for free groups.
A negative answer would imply a negative answer to M. Gromov's Question, Question 1.1 of \cite{bestvinaproblems}, regarding
necessary and sufficient conditions for a group to be word hyperbolic.
See Section 8 of \cite{klsurvey} for a
discussion of the connection between Gromov's question and convex cocompactness.

Question 6 of \cite{kentleininger} is a natural test question for Question 1.5 of \cite{FMcc} as the necessary and
sufficient condition for an element in $\pi_1(S,z)$ to be pseudo-Anosov as an element of $\Mod(S,z)$ is a topological
one, and not \textit{a priori} related to any algebraic structure.
I. Kra discovered this necessary and sufficient condition \cite{Kra}---see Theorem
\ref{T:Kra} here---and his proof of sufficiency is Teichm\"uller theoretic (necessity is
obvious).
We give an alternative proof here
based entirely on topological and group theoretic considerations, see \S \ref{S:stabilizer section}.

The class of groups covered by Theorem \ref{T:convex ccpt 1} also includes the test case proposed by Farb and Mosher in
Question 1.6 of \cite{FMcc} (see also Problem 3.5 of \cite{Mosherproblems}).  These are the finitely generated
subgroups of K. Whittlesey's groups.
Recall that Whittlesey's groups are normal purely pseudo-Anosov subgroups of the mapping class groups of the sphere with $n \geq 5$ punctures and of the closed genus--$2$ surface.

\begin{corollary*} \label{C:Whittlesey}
Whittlesey's groups are locally convex cocompact: finitely generated subgroups are convex cocompact.
\end{corollary*}
\begin{proof}[Proof of the Corollary from Theorem \ref{T:convex ccpt 1}]
It suffices to prove the theorem for Whittlesey's subgroups of $\Mod(S_{0,n})$,  as there is a surjection
\[
\Mod(S_{2,0}) \to \Mod(S_{0,6})
\]
with order--two central kernel \cite{birmanhilden},
 and an isometry of Teichm\"uller spaces which is equivariant with respect to this ``virtual isomorphism.''

We can view any one of the punctures of $S_{0,n}$ as being obtained from $S_{0,n-1}$ by removing a marked point $z$.
There are thus $n$ different Birman sequences and so $n$ surjective homomorphisms $\Mod(S_{0,n}) \to \Mod(S_{0,n-1})$.
The intersection of these kernels is Whittlesey's group, and hence lies in $\pi_1(S_{0,n-1}) < \Mod(S_{0,n})$.
Any finitely generated subgroup of Whittlesey's group is thus also a finitely generated purely pseudo-Anosov subgroup
of $\pi_1(S_{0,n-1})$. Since $n \geq 5$, Theorem \ref{T:convex ccpt 1} implies that such a subgroup is convex cocompact.
\end{proof}

\bigskip

The proof of Theorem \ref{T:convex ccpt 1} relies on a characterization of convex cocompactness discovered by the first two authors \cite{kentleininger} and, independently, by
Hamenst\"adt \cite{hamenstadt} in terms of the action
on the \textit{curve complex} $\C$.
We are thus lead to a study of the action of $\pi_1(S,z) < \Mod(S,z)$ on
the curve complex of $(S,z)$ (equivalently, that of $S \setminus \{z \}$), and there is an interesting observation
regarding this action that we now describe.

If $S$ is closed ($m = 0$), then there is a map of curve complexes $\Pi \co \C(S,z) \to \C(S)$; see \S \ref{S:curve
complex}.
In general, the desired map is not globally well-defined as essential simple closed curves on $(S,z)$ may become
peripheral in $S$.  In this case we restrict $\Pi$ to the largest subcomplex on which it is well-defined, denoted $\widehat
\C(S,z)$. The fibers of this map are invariant under the action of $\pi_1(S,z) < \Mod(S,z)$, and have a simple
geometric description:

\begin{theorem_mt1*}
The fiber of $\Pi$ over a point in the interior of a simplex $v \subset \C(S)$ is $\pi_1(S,z)$--equivariantly
homeomorphic to the tree $T_v$ determined by $v$.
\end{theorem_mt1*}

The tree $T_v$ is the tree dual to the multi-curve $v$ equipped with its action by $\pi_1(S,z)$; see \S \ref{S:tree
section}.  This is the Bass--Serre tree for the splitting of $\pi_1(S,z)$ determined by the multi-curve
$v$.

A consequence of the theorem is the following fact due to J. Harer \cite{Harervcd}---a new proof of this is due to
A. Hatcher and K. Vogtmann \cite{hatchervogtmann}.

\begin{corollary} [Harer]
With the polyhedral topologies, $\C(S)$ and $\widehat \C(S,z)$ are homotopy equivalent.
\end{corollary}

The corollary is proven using a section of $\Pi$ described by Harer and performing the straight--line homotopy  to
the section along the fibers given by Theorem \ref{T:main tree 1}.  See \S \ref{S:tree section}.

\bigskip

\noindent {\bf Acknowledgments.}  The authors thank Yair Minsky, Ursula Hamenst\"adt, Alan Reid, and Ben Wieland for
helpful and interesting conversations.  The second author thanks the Max-Plank-Institut f\"ur Mathematik in Bonn for
its hospitality during part of this work.   We also thank Andy Putman for his careful reading and
comments on an earlier version of this paper and the referee for providing helpful organizational advice.


\section{Definitions and conventions} \label{S:background}

We have chosen to work with the surface $(S,z)$ marked with $z$ rather than $S \setminus \{z\}$ as this is generally
more convenient.
However, occasionally our arguments are clarified by working with $S \setminus \{ z \}$.
When this is the case, we refer to the puncture obtained by removing $z$ as the \textbf{$z$--puncture}.

We say that a closed curve in $S$ is \textbf{nontrivial} if it is homotopically nontrivial in $S$ and
\textbf{essential} if it is nontrivial and nonperipheral, that is, not homotopic into every neighborhood of a puncture.
These definitions are extended to $(S,z)$ by defining a closed curve in $(S,z)$ to be a closed curve in $S$ which is
contained in $S \setminus \{z\}$.  A closed curve in $(S,z)$ is then \textbf{nontrivial} (respectively, \textbf{essential}) if it is so in $S \setminus
\{z\}$. Isotopy in $(S,z)$ means isotopy in $S$ fixing $z$.
Thus nontrivial and essential simple closed curves in
$(S,z)$ are isotopic if and only if they are isotopic in $S \setminus \{z\}$.

We fix a complete finite area hyperbolic metric on $S$ and let $p \co \widetilde S \to S$ denote the universal covering.
The hyperbolic metric on $S$ pulls back to one on $\widetilde S$ making $\widetilde S$ isometric to the hyperbolic plane.

We view $\pi_1(S)$ as the group of covering transformations of the universal covering $p\co \widetilde S \to S$ and fix
this action once and for all.  A point $\widetilde z \in p^{-1}(z)$ determines an isomorphism of $\pi_1(S)$ with the
fundamental group $\pi_1(S,z)$.
We fix a basepoint $\widetilde z \in p^{-1}(z)$, and hence an
isomorphism $\pi_1(S) \cong \pi_1(S,z)$.

\subsection{Mapping class groups and Birman's sequence}

The \textbf{mapping class group} of $S$ is the group
$\Mod(S) = \pi_0(\Diff^+(S))$, where $\Diff^+(S)$ is
the group of orientation preserving diffeomorphisms of $S$ that fix each of the punctures. We define $\Mod(S,z)$ to be
$\pi_0(\Diff^+(S,z))$, where $\Diff^+(S,z)$ is the group of orientation preserving diffeomorphisms of $S$ that fix
each puncture and that also fix $z$.  There is a canonical isomorphism $\Mod(S,z) \cong \Mod(S \setminus \{z\})$.

Birman's exact sequence
\cite{birmansequence,birmanbook} relates the mapping class group of $S$ with that of $(S,z)$ and $\pi_1(S,z)$. Namely
\[
1 \to \pi_1(S,z) \to \Mod(S,z) \to \Mod(S) \to 1.
\]

To describe the inclusion $\pi_1(S,z) \to \Mod(S,z)$ concretely, we first represent an element of $\pi_1(S,z)$ by a
loop $\gamma$ based at $z$. Writing $\gamma \co [0,1] \to S$ with $\gamma(0) = \gamma(1) = z$, there is an isotopy $h_t \co S
\to S$ such that $h_0 = \Id_S$ and $\gamma(t) = h_t(z)$ for all $t \in [0,1]$.
Since $h_1(z) = z$, the map $h_1$ determines a
mapping class in $\Mod(S,z)$, and this is the image of $\gamma$ in $\Mod(S,z)$ in the exact sequence.
It is clear that the isotopy $h_t$ may be constructed so that the diffeomorphism $h_1$ is supported on any given neighborhood of the curve $\gamma \subset S$.

For clarity, we write $h_\gamma$ for the diffeomorphism or mapping class associated to $\gamma \in \pi_1(S,z)$.

Any element $f \in \Mod(S)$ determines an outer automorphism $f_*$ of $\pi_1(S,z)$,
and this defines a homomorphism $\Mod(S) \to \Out(\pi_1(S,z))$---the codomain is $\Out(\pi_1(S,z))$ rather than $\Aut(\pi_1(S,z))$ as representatives of $f$ need not fix the basepoint.
The Dehn--Nielsen Theorem states that this is an isomorphism onto an index--two
subgroup when $S$ is closed (see Stillwell's appendix to \cite{dehnstillwell}), but in general it is only injective \cite{stallingsPV}.
Working with $\Mod(S,z)$ erases the difficulty of
moving basepoints and there is an injection $\Mod(S,z) \to \Aut(\pi_1(S,z))$.  Since $\pi_1(S,z)$ has trivial center, $\pi_1(S,z) \cong \Inn(\pi_1(S,z))$ and Birman's exact sequence injects into the classical short exact sequence associated to $\Aut(\pi_1(S,z))$:
\[
\xymatrix{
1 \ar[r] & \pi_1(S,z) \ar[r]\ar[d] & \Mod(S,z) \ar[r]\ar[d] & \Mod(S) \ar[r]\ar[d] & 1\\
1 \ar[r] & \Inn(\pi_1(S,z)) \ar[r] & \Aut(\pi_1(S,z)) \ar[r] & \Out(\pi_1(S,z)) \ar[r] & 1\\}
\]
All of the squares commute, and, for any $\alpha,\gamma \in \pi_1(S,z)$, the first one becomes
\begin{equation} \label{E:mod to inner}
(h_\gamma)_*(\alpha) = \gamma \alpha \gamma^{-1} .
\end{equation}

\subsection{Curve complexes} \label{S:curve complex}

The \textbf{simplicial curve complex of $S$} will be denoted $\C^\Delta(S)$, the simplicial curve complex of $(S,z)$ similarly denoted $\C^\Delta(S,z)$.  We view these abstract
simplicial complexes as partially ordered sets of simplices, where a $k$-simplex is determined by a set $v = \{
v_0,\ldots,v_k \}$ of $k+1$ distinct isotopy classes of pairwise disjoint essential simple closed curves.

Following the nomenclature of the subject, we let $\C(S)$ and $\C(S,z)$ denote the \textbf{curve complexes of $S$} and \textbf{$(S,z)$}, respectively.  These are geodesic metric spaces obtained by isometrically gluing regular Euclidean simplices with all edge lengths equal to one according to the combinatorics of the associated abstract simplicial complex (compare \cite[I.7]{BH}).  The necessity for distinguishing between $\C$ and $\C^\Delta$ will soon become apparent.

In the case that $\xi(S) = 1$, $\C^\Delta(S)$ is zero dimensional and one often makes a separate definition for the curve complexes in these cases.
However, \textit{we do not do this here}, and so considering $\C(S)$ a geodesic metric space is nonsensical---a geodesic metric
would have to assign an ``infinite distance'' to any two points.  In this special case, we simply treat $\C(S)$
as a countable set of points.
As we will only consider metric properties of $\C(S,z)$, and not $\C(S)$,
this is not a serious issue.  For this reason, and so we need not continue to comment on the special case, we
discuss the relevant metric geometry of $\C(S,z)$ only.

The $1$-skeleton $\C^1(S,z)$ is itself a metric space (with the induced geodesic metric), and the inclusion into
$\C(S,z)$ is a quasiisometry. Because a geodesic in $\C^1(S,z)$ between vertices has a combinatorial description
as a sequence of adjacent vertices, we may mix combinatorial and geometric arguments in the metric space
$\C^1(S,z)$. We will therefore work with the metric on $\C^0(S,z)$ induced by the inclusion into $\C^1(S,z)$, which
takes integer values only.

We will confuse notation as follows:  a simplex in $\C^\Delta$ will be considered simultaneously a subset
of $\C$ as well as a union of curves on the surface (further confusing an isotopy class with a representative).  Thus we may write $v \in \C^\Delta(S)$, $v \subset \C(S)$, or $v \subset
S$, the context determining the particular meaning.

For clarity, we will typically denote simplices of $\C^\Delta(S,z)$ by $u = \{u_0,\ldots,u_k\}$ and simplices of
$\C^\Delta(S)$ by $v = \{ v_0,\ldots,v_k \}$. If $v$ is a simplex in $\C^\Delta(S)$, we will write $[v]$ for the
geodesic representative of $v$, which is a union of pairwise disjoint embedded simple closed geodesics.

\subsection{Forgetful projection} \label{S:projection}

Any simple closed curve $u$ in $(S,z)$ can be viewed as a curve in $S$ which we denote $\Pi(u)$.  However, if $S$ has punctures then an essential curve in $(S,z)$ may become
inessential in $S$ by becoming peripheral. The only time this can happen is when the curve is the boundary of a
once--punctured disk $(Y,z) \subset (S,z)$ containing the marked point.
We call such a curve $u$ in $(S,z)$
\textbf{preperipheral}.
A simplex $u = \{u_0,\ldots,u_k\} \in \C^\Delta(S,z)$ is
called \textbf{preperipheral} if one of its vertices $u_i$ is a preperipheral curve in $(S,z)$, and \textbf{nonpreperipheral}
otherwise.

We define the \textbf{nonpreperipheral subcomplex} of $\C^\Delta(S,z)$ by
\[
\widehat C^\Delta(S,z) = \{ u \in \C^\Delta(S,z) \, | \, u \mbox{ is nonpreperipheral} \}.
\]
There is now a well-defined simplicial map
\[
\Pi_\Delta \co  \widehat \C^\Delta(S,z) \to \C^\Delta(S)
\]
determined by forgetting the marked point $z$.  We also write
\[
\Pi \co \widehat \C(S,z) \to \C(S)
\]
for the induced map on metric spaces (or sets if $\xi(S) = 1$).

\begin{lemma} \label{L:small collapse}
If $u \subset \widehat \C(S,z)$ is a $k$--simplex, then $\Pi(u)$ has dimension at least $k-1$.

If $\Pi|_u$ is noninjective (so that $\Pi(u)$ has dimension $k-1$) and after reordering the vertices of $u$ we have $\Pi(u_0)=\Pi(u_1)$, then the curves $u_0$ and $u_1$ cobound an
annulus $(Y,z) \subset (S,z)$ containing $z$.
\end{lemma}

\begin{proof}
Every simplex $u \in \widehat C^\Delta(S,z)$ is contained in a maximal simplex $u'$ in $C^\Delta(S,z)$, which
determines a pants decomposition (of $S \setminus \{z\}$).
So $u'$ has dimension $\xi(S \setminus \{z \}) - 1=
\xi(S)$.
Since $\xi(S) \geq 1$, an easy argument allows us to assume that $u' \in \widehat \C(S,z)$.

Every component of $(S \setminus \{ z \}) \setminus u'$ is a pair of pants, exactly one of which contains the
$z$--puncture.
So the corresponding component $(Y,z) \subset (S,z)$ of $S \setminus u'$ is an annulus with two boundary
components.  After reordering the vertices we may assume these are $u_0,u_1 \in u'$.  Clearly $\Pi(u')$ is a
pants decomposition of $S$, and so has dimension $\xi(S) -1$, one less than that of $u$.  Since $\Pi|_{\{u_0,u_1\}}$
is not injective, it follows that $\Pi|_{u}$ will be noninjective if and only if it contains both $u_0$ and $u_1$,
which thus cobound the annulus $(Y,z)$.  In any case, the dimension can be at most one less than that of $u$.
\end{proof}

We say that a simplex $u$ of $\widehat \C(S,z)$ is {\bf injective} if $\Pi|_u$ is injective and {\bf noninjective}
otherwise.

\section{Subsurfaces} \label{S:subsurfaces}

Consider a compact $\pi_1$--injective subsurface $(Y,z) \subset (S,z)$ with $Y \not \cong D^2$ and $\pi_1(Y,z) <
\pi_1(S,z)$ a proper subgroup. The boundary $\partial Y$ is a disjoint union of nontrivial simple closed curves in
$(S,z)$.  It can happen that some of the components in $\partial Y$ are isotopic to each other in $(S,z)$ and that some of
the components are peripheral in $(S,z)$.  We let $\partial_0 Y$ denote the simplex in $\C^\Delta(S,z)$ obtained by
identifying pairs of components of $\partial Y$ that are isotopic in $(S,z)$ and forgetting the peripheral components.

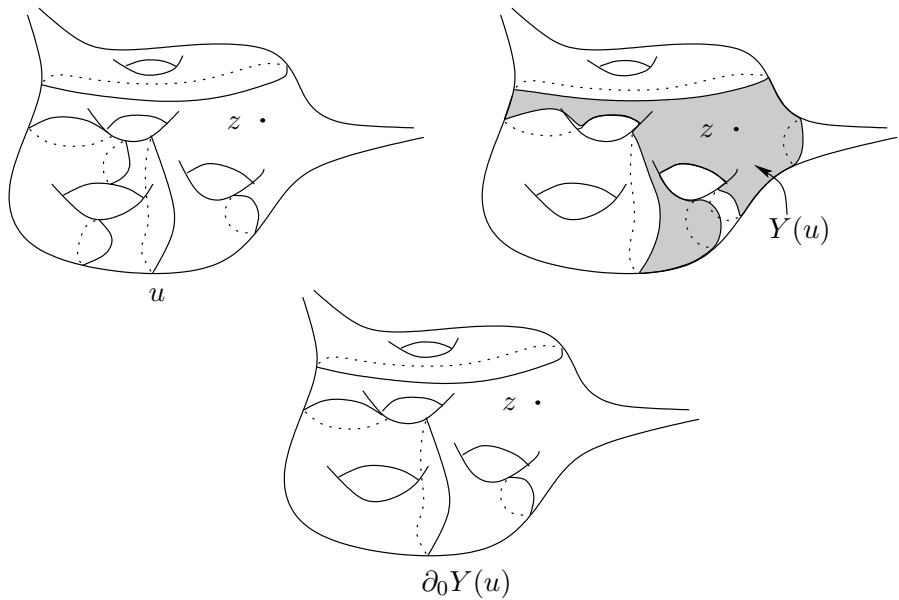
\begin{figure}[h!]
\begin{center}
\input{subsurface.pstex_t}
\end{center}
\caption{A simplex $u \in \C^\Delta(S,z)$, its subsurface $(Y(u),z)$, and $\partial_0 Y(u) \subseteq u$.}
\label{F:subsurface}
\end{figure}

Let $u \in \C^\Delta(S,z)$ be any simplex.
We construct a subsurface $(Y(u),z)
\subset (S,z)$ as follows.  Remove from $S$ a small open tubular neighborhood of $u$ and small open cusp neighborhoods
of the punctures (in particular, these neighborhoods should be pairwise disjoint and not contain $z$), and let $Y(u)$ be the component of
this subsurface containing $z$.  Note that $\partial_0 Y(u) \subseteq u$ is a potentially proper face---see Figure
\ref{F:subsurface} for an example.  We can iterate this process of taking subsurface, then $\partial_0$, but it
immediately stabilizes: $Y( \partial_0 Y(u))$ is isotopic fixing $z$ to $Y(u)$.

\begin{lemma} \label{L:subsurfacesimplextype}
If $u \in \C(S,z)$, then $(Y(u),z)$ is an annulus (containing $z$) if and only if $u$ is either preperipheral or nonpreperipheral and noninjective.
Furthermore, in this case $u$ is preperipheral if and only if the boundary components of $Y(u)$ are peripheral in $S$.
\end{lemma}
\begin{proof}
If $u \in \widehat \C^\Delta(S,z)$, the lemma follows immediately from Lemma 2.1 and the definitions.

All that remains is to prove that for a preperipheral simplex $u \in \C^\Delta(S,z)$, $Y(u)$ is an annulus for which the boundary components are peripheral in $S$.
Let $u_0 \in u$ be the preperipheral vertex.  It follows that $u$ bounds a once-punctured disk containing $z$, that is, $Y(u_0)$ is an annulus of the required type.  Now let $u_1 \in u$ be any other vertex.  The curve $u_1$ cannot be contained in $(Y(u),z)$, since $Y(u)-\{z\}$ is a pair of pants and thus contains no essential simple closed curves.  Therefore, all other components of $u$ lie outside $Y(u_0)$, and so $Y(u_0) = Y(u)$, completing the proof.
\end{proof}

\subsection{From simplices to groups} \label{S:simplex to group}

To each subsurface $(Y,z) \subset (S,z)$ described in the previous section we can associate its fundamental group $\pi_1(Y,z) < \pi_1(S,z)$.
We let $\D$ denote the collection of all such subgroups and define a map
\[
\Gamma \co  \C^\Delta(S,z) \to \D
\]
by declaring that $\Gamma(u) = \pi_1(Y(u),z)< \pi_1(S,z)$.  Note that while $(Y(u),z)$ is only defined up to isotopy fixing $z$
(since $u$ is), $\Gamma(u)$ is a well defined subgroup of $\pi_1(S,z)$.

The set $\D$ admits a natural partial order by inclusion as well as an action of $\pi_1(S,z)$ by conjugation.
\begin{proposition} \label{P:complex to groups}
The map $\Gamma$ is a $\pi_1(S,z)$--equivariant order--reversing surjection.
\end{proposition}
\begin{proof}
Given a subgroup $\pi_1(Y,z) < \pi_1(S,z)$ in $\D$, we have $\pi_1(Y(\partial_0 Y),z) = \pi_1(Y,z)$, and hence $\Gamma(\partial_0 Y) = \pi_1(Y,z)$.
So $\Gamma$ is
surjective. If $u \subseteq u'$ is a face, then $Y(u) \supseteq Y(u')$, and so $\pi_1(Y(u),z) \geq \pi_1(Y(u'),z)$.
In other words, $\Gamma(u) \geq \Gamma(u')$.
So $\Gamma$ is order--reversing.

If $\gamma \in \pi_1(S,z)$, then $Y(h_\gamma u) = h_\gamma(Y(u))$.
We have
\begin{align*}
\Gamma(h_\gamma u)
= \pi_1(Y(h_\gamma u), z)
& = \pi_1(h_\gamma (Y(u)),z) \\
& = \gamma \pi_1(Y(u),z) \gamma^{-1} \quad \quad \mathrm{by\ \eqref{E:mod to inner}}\\
&= \gamma \Gamma(u) \gamma^{-1}
\end{align*}
and it follows that $\Gamma$ is equivariant.
\end{proof}

\subsection{From groups to convex hulls} \label{S:group to hull}

In the following, all stabilizers are taken with respect to actions of $\pi_1(S)$ and are denoted $\Stab(\cdot)$.
Our choice of isomorphism $\pi_1(S) \cong \pi_1(S,z)$ allows us to view $\D$ as a collection of subgroups of $\pi_1(S)$
acting on the hyperbolic plane $\widetilde S$.
The type of group $\Gamma(u)$ is reflected by the type of
the simplex $u$.  More precisely, from the definitions and Lemma \ref{L:subsurfacesimplextype} we have the
following:
\begin{enumerate}
\item[(i)] $\Gamma(u)$ is a nonabelian free group if and only if $u$ is nonpreperipheral and $u$ is injective.
\item[(ii)] $\Gamma(u)$ is cyclic generated by a hyperbolic element if and only if $u$ is nonpreperipheral and $u$ is noninjective.
\item[(iii)] $\Gamma(u)$ is cyclic generated by a parabolic element if and only if $u$ is preperipheral.
\end{enumerate}

In cases (i) and (ii), we let $\hull(\Gamma(u))$ denote the convex hull of the limit set $\Lambda_{\Gamma(u)}$ of $\Gamma(u)$ acting on
$\widetilde S \cong \mathbb H^2$.  In case (ii), $\Gamma(u)$ is cyclic and $\hull(\Gamma(u))$ is the axis of the
elements of $\Gamma(u)$.  Alternatively, $\hull(\Gamma(u))$ is a component of $p^{-1}([\Pi(u)])$ (recall that $[\Pi(u)]$ is the geodesic representative of $\Pi(u)$).
Note that $\Gamma(u) = \Stab(\hull(\Gamma(u)))$ since the natural map of the quotient $\hull(\Gamma(u))/\Gamma(u) = [\Pi(u)]$ into
$S$ is injective.

If $\Gamma(u)$ is not cyclic (case (i)), then $\hull(\Gamma(u))$ is a subsurface of $\widetilde S$ bounded by geodesics.   Consider the
isotopy (not fixing $z$) from $\Pi(u)$ to the geodesic representative $[\Pi(u)]$ in $S$.  This takes $Y(u)$ to a compact
core for the quotient of the interior $\intr(\hull(\Gamma(u)))/\Gamma(u) \subset S$.  This can be seen by lifting the
isotopy to a $\pi_1(S)$--equivariant isotopy in $\widetilde S$ taking $p^{-1}(\Pi(u))$ to $p^{-1}([\Pi(u)])$.  Because
$\intr(\hull(\Gamma(u)))/\Gamma(u)$ injects into $S$ and since $\Stab(\hull(\Gamma(u))) =
\Stab(\intr(\hull(\Gamma(u))))$ it follows that $\Gamma(u) = \Stab(\hull(\Gamma(u)))$.

In case (iii), the convex hull of the limit set of $\Gamma(u)$ is empty since the limit set is a single point.
Here we define $\hull(\Gamma(u))$ to be a $\Gamma(u)$--invariant horoball, chosen as follows.
We choose a $\pi_1(S)$--invariant family of horoballs, one centered at each parabolic fixed point, with the property that for
any horoball in the family
\begin{itemize}
\item[(I)\ \ ] the quotient by the stabilizer embeds as a cusp neighborhood of the associated puncture in $S$

\item[ \textit{and}]
\item[(II)\ ] the boundary of the cusp neighborhood (which is the quotient of the boundary horocycle by the
stabilizer) is disjoint from every simple closed geodesic in $S$.
\end{itemize}
That this is possible is a well--known consequence
of J{\o}rgensen's inequality (more precisely, the Shimizu--Leutbecher Lemma \cite[II.C]{Maskit}).  In Section \ref{S:convex cocompactness} we impose
tighter restrictions on these horoballs, but, for now, this suffices.

Note that we also have $\Stab(\hull(\Gamma(u))) = \Gamma(u)$ in case (iii) since $\Gamma(u)$ is necessarily a
\textit{maximal} parabolic subgroup.
We therefore have
\begin{lemma} \label{L:stabilizer of hull}
For every $\Gamma(u) \in \D$, $\Stab(\hull(\Gamma(u))) = \Gamma(u)$. \qed
\end{lemma}

Let $N(\Gamma(u))$ denote the normalizer of $\Gamma(u)$ in $\pi_1(S)$.  Lemma \ref{L:stabilizer of hull} yields
\begin{proposition} \label{P:self normalize}
For every $\Gamma(u) \in \D$, $N(\Gamma(u)) = \Gamma(u)$.
\end{proposition}
\begin{proof}
As is true for any subgroup of $\pi_1(S)$, we have
\[
\Gamma(u) \leq N(\Gamma(u)) \leq \Stab(\Lambda_{\Gamma(u)})
\]
and
\[
\Stab(\Lambda_{\Gamma(u)}) = \Stab(\hull(\Gamma(u))).
\]
Combining these relations with Lemma \ref{L:stabilizer of hull} we have
\[
\Gamma(u) \leq N(\Gamma(u)) \leq \Gamma(u),
\]
and the proposition follows.
\end{proof}

We have associated to each  $\Gamma(u) \in \D$ a convex set $\hull(\Gamma(u))$ and we let $\H$ denote the set of all
such convex sets.  We order $\H$ by inclusion and note that $\pi_1(S)$ clearly acts on $\H$ as  convex hulls of limit sets are natural and our set of
horoballs is $\pi_1(S)$--invariant.
\begin{lemma} \label{L:groups to hulls}
The map $\hull \co \D \to \H$ is a $\pi_1(S)$--equivariant order--preserving surjection.
\end{lemma}
\begin{proof}
First observe that $\Gamma(u) \leq \Gamma(u')$ implies $\Lambda_{\Gamma(u)} \subseteq \Lambda_{\Gamma(u')}$.

If neither group $\Gamma(u)$ nor $\Gamma(u')$ is parabolic, then it immediately follows that $\hull(\Gamma(u)) \subseteq \hull(\Gamma(u'))$. If both
groups are parabolic, then we must have $\Gamma(u) = \Gamma(u')$, for the parabolic subgroups of $\D$ are maximal
parabolic subgroups, and so $\hull(\Gamma(u)) = \hull(\Gamma(u'))$.

The only remaining case is when $\Gamma(u)$ is parabolic and $\Gamma(u')$ is not.  Since $\Gamma(u) < \Gamma(u')$, some $\Gamma(u)$--invariant horoball is contained in $\hull(\Gamma(u'))$.  Since $\hull(\Gamma(u'))$ is
bounded by geodesics that descend to simple closed geodesics in $S$, and these geodesics are disjoint from the
boundaries of the cusps $\hull(\Gamma(u))/\Gamma(u)$ by construction, it follows that $\hull(\Gamma(u)) \subseteq
\hull(\Gamma(u'))$.

The map is surjective by construction.
\end{proof}

Proposition \ref{P:complex to groups} and Lemma \ref{L:groups to hulls} have the following corollary.
\begin{corollary} \label{C:complex to convex sets}
\[
\hull \circ \Gamma \co  \C^\Delta(S,z) \to \H
\]
is a $\pi_1(S)$--equivariant order--reversing surjection. \qed
\end{corollary}

\section{Stabilizers and Kra's Theorem} \label{S:stabilizer section}

\begin{theorem} \label{T:simplex stab}
If $u \in \C^\Delta(S,z)$ is any simplex, then
\[
\Stab(u) = \Gamma(u)
\]
with respect to the action determined by the inclusion $\pi_1(S,z) < \Mod(S,z)$.
\end{theorem}
We note that the stabilizer in $\pi_1(S,z) < \Mod(S,z)$ of a simplex $u$ fixes that simplex pointwise
since $\pi_1(S,z)$ acts trivially on homology, and so by Theorem 1.2 of \cite{ivanov}, it consists entirely of pure
mapping classes.
\begin{proof}
One inclusion is obvious.  Namely,
\[
\Gamma(u) = \pi_1(Y(u),z) \leq \Stab(u)
\]
for if $\gamma \in \pi_1(Y(u),z)$, then the diffeomorphism $h_\gamma$ can be chosen to be supported on $Y(u)$.  In particular, $h_\gamma$ fixes $u$.

To prove the other inclusion, observe that the $\pi_1(S,z)$--equivariance of
$\Gamma$ implies that the stabilizer of $u$ is contained in the stabilizer of $\Gamma(u)$.   Since
$\pi_1(S,z)$ acts by conjugation on $\D$ this means that
\[
\Stab(u) \leq N(\Gamma(u)).
\]
Proposition
\ref{P:self normalize} implies $N(\Gamma(u)) = \Gamma(u)$ and so
\[
\Stab(u) \leq \Gamma(u)
\]
as required.

\end{proof}

We obtain Kra's Theorem as a corollary:
\begin{theorem}[Kra] \label{T:Kra}
Given $\gamma \in \pi_1(S,z)$, $h_\gamma \in \Mod(S,z)$ is pseudo-Anosov if and only if $\gamma$ is filling.
\end{theorem}
An element of $\pi_1(S,z)$ is filling if every representative loop nontrivially intersects every essential closed curve
in $S$.  We note that every representative loop of an element $\gamma \in \pi_1(S,z)$ intersects every essential closed
curve in $S$ if every closed curve in the \textit{free} homotopy class of $\gamma$ does.
\begin{proof}
Since $h_\gamma$ is pure, it is pseudo-Anosov if and only if it does not stabilize any simplex. By Theorem
\ref{T:simplex stab} this happens if and only if $\gamma$ is not in any group $\pi_1(Y(u),z)$ for any $u \in
\C^\Delta(S,z)$.

We claim that $\gamma$ is not in any group $\pi_1(Y(u),z)$ for any $u \in C^\Delta(S,z)$ if and only if $\gamma$ is
filling.  To see this, first observe that $\gamma$ is in $\pi_1(Y(u),z)$ if and only if it can be realized disjoint
from $\partial_0 Y(u)$.  If $\partial_0 Y(u)$ is not preperipheral, then $\gamma$ has a representative disjoint from
any component of $\Pi(\partial_0 Y(u))$, which is an essential curve in $S$.  If $\partial_0 Y(u)$ is preperipheral,
then $\gamma$ is a peripheral loop, so has a representative disjoint from a representative of any essential closed
curve in $S$.
\end{proof}

\section{Purely pseudo-Anosov subgroups} \label{S:purely pa}

Let $G < \pi_1(S,z)$ be a finitely generated subgroup, purely pseudo-Anosov when considered a subgroup of $\Mod(S,z)$.
We view $G$ as a subgroup of $\pi_1(S)$ acting as a fuchsian group on the hyperbolic plane $\widetilde S$ and make uniform estimates  to be used in the
proof of Theorem \ref{T:convex ccpt 1}.

By Theorem \ref{T:Kra}, the free homotopy class in $S$ defined by any nontrivial element of $G$ fills $S$. Let $\Sigma =
 \hull(G)/G$ be the quotient hyperbolic surface with geodesic boundary and let
\[
p_0\co \hull(G) \to \Sigma
\]
be the covering projection.

The surface $\Sigma$ is compact. To see this, note that
finite generation of fuchsian groups is equivalent to geometric finiteness (see Theorem 10.1.2 of \cite{beardon}) and 
that $G$ is purely hyperbolic as a fuchsian group acting on $\widetilde S$, as every element corresponds to a filling loop on $S$.

The inclusion $\hull(G) \to \widetilde S$ induces an immersion $f \co \Sigma \to S$ with
$f_*(\pi_1(\Sigma)) = G$.  We collect our maps into a commuting diagram:
\[
\xymatrix{
\widetilde{\Sigma} \ar[r] \ar[d]_{p_0} & \widetilde{S} \ar[d]^p\\
\Sigma \ar[r]^f & S\\}
\]

Since every nontrivial conjugacy class in $G$ is filling, for any geodesic $[v]$ in $S$, $f^{-1}([v])$ cuts $\Sigma$ into disks.
To see this, note that if this were not the case, then there would be a nontrivial loop $\gamma \in \pi_1(\Sigma)$
disjoint from $f^{-1}([v])$, and so $f_*(\gamma)$ could not be filling as it would be disjoint from $v$. Moreover, as we
will see, the family of arcs of $f^{-1}([v])$ {\em as $v$ ranges over all of $\C^0(S)$} is a precompact family.


It will be convenient in the proof of the lower bound for Theorem \ref{T:quasiisometric embedding} to demonstrate
precompactness of the family in a slightly larger surface.  Namely, let $\Sigma_1 = N_1(\hull(G))/G$ denote the
quotient of the $1$-neighborhood of $\hull(G)$ by $G$. This adds a collar of width one to each boundary
component of $\Sigma$. There is an obvious extension of $f$ to $\Sigma_1 \supset \Sigma$ that we still denote
$f \co \Sigma_1 \to S$. Let ${\mathcal A}$ denote the set of all arcs of $f^{-1}([v])$ in $\Sigma_1$ as $v$ ranges over all
of $\C^0(S)$.

\begin{proposition} \label{P:arcs compact}
There are only finitely many isotopy classes in ${\mathcal A}$ and there is a uniform upper bound on the length of any arc in ${\mathcal A}$.
\end{proposition}

\begin{proof}
The first statement clearly follows from the second.

Suppose that there is a sequence $\{ v_n \} \subset \C^0(S)$ and components $L_n \subset
f^{-1}([v_n])$ so that $\ell(L_n) \to \infty$.
After passing to a subsequence, we may assume
that $[v_n]$ has a Hausdorff limit $\lambda$, a geodesic lamination on $S$.  Let $\lambda'$ be the maximal
measurable sublamination of $\lambda$, obtained from $\lambda$ by throwing away all non-closed
isolated leaves.

After passing to a further subsequence if necessary, we assume that $L_n$ has a Hausdorff limit, necessarily contained in $f^{-1}(\lambda)$.
Since $\ell(L_n) \to \infty$, there is a connected
geodesic lamination $\kappa$ contained in this limit and $f(\kappa)$ is a component of $\lambda'$---the components of
$\lambda'$ are exactly the sublaminations.
 If $\kappa$ is a simple closed geodesic, then $f(\kappa) \subset
\lambda'$ represents a conjugacy class in $G$ which is not filling and we obtain a contradiction.
So $\kappa$ is not a simple closed
geodesic.

Let $W_\kappa$ denote the supporting subsurface of $\kappa$---the smallest open, locally convex subsurface containing
$\kappa$---and let $Y_{f(\kappa)}$ be the supporting subsurface of $f(\kappa)$.  Since the supporting subsurface is
determined by the preimage of the lamination in the universal cover, it easily follows that
\[
f(W_\kappa) = Y_{f(\kappa)}.
\]
Thus the surface $f(W_\kappa)$ is the component of the supporting subsurface of $\lambda'$ containing $f(\kappa)$ and
$f(\partial W_\kappa)$ is disjoint from $\lambda'$.  This is impossible since every curve representing a conjugacy
class of $G$ intersects every lamination in $S$, and any component of $f(\partial W_\kappa)$ represents a conjugacy
class in $G$.
\end{proof}

Note that each component of $\Sigma_1 \setminus f^{-1}([v])$ is a not only a disk, but a disk with uniformly bounded
diameter (with respect to the induced path metric): it is convex, and has uniformly bounded circumference.  This
implies the same statement for the disks $\overline{N_1(\hull(G)) \setminus p^{-1}([v])}$, which are just the
disks of $\overline{p_0^{-1}(\Sigma_1 \setminus f^{-1}([v]))}$. Moreover, since $N_1(\hull(G)) \cap
p^{-1}([v]) = p_0^{-1}(\Sigma_1 \cap f^{-1}([v]))$, the lemma tells us that these arcs have uniformly bounded
diameter also.

The components of $p^{-1}([v])$ and the closures of the components of $\widetilde S \setminus p^{-1}([v])$ are
precisely the convex hulls $\hull(\Gamma(u))$ for $\Gamma(u) \in \D$ with $\Pi_\Delta(u) = v$.  With the remarks of the previous paragraph and
Proposition \ref{P:arcs compact}, this implies the following.

\begin{corollary} \label{C:uniform diameter}
There exists a $D > 0$ so that for any simplex $u \in \widehat \C^\Delta(S,z)$,
\[
\diam(\hull(\Gamma(u)) \cap N_1(\hull(G))) \leq D.
\]
\end{corollary}
\vspace{-19pt}
\qed

\section{Convex cocompactness} \label{S:convex cocompactness}

\begin{theorem} \label{T:convex ccpt 1}
If $\xi(S) \geq 1$ and $G < \pi_1(S,z)$ is finitely generated and purely pseudo-Anosov as a subgroup of $\Mod(S,z)$,
then $G$ is convex cocompact.
\end{theorem}

Fix a purely pseudo-Anosov subgroup $G < \pi_1(S,z) < \Mod(S,z)$ and let $\hull(G) \subset N_1(\hull(G)) \subset \widetilde S$ and $f  \co  \Sigma_1 \to S$ be as in the previous section.  We assume that for any preperipheral simplex $u \in \C^\Delta(S,z)$ we have chosen $\hull(\Gamma(u))$ so that
\[
\hull(\Gamma(u)) \cap N_1(\hull(G)) = \emptyset.
\]
This is possible since $f(\Sigma_1)$ is a compact subset of $S$, and there are only finitely many such conjugacy
classes of subgroups $\Gamma(u)$ (recall from \S \ref{S:group to hull} that the $\Gamma(u)$ are precisely the maximal parabolic subgroups and the $\hull(\Gamma(u))$ are invariant horoballs).

From this choice and Corollary \ref{C:uniform diameter}, we obtain the following refinement of that corollary.

\begin{proposition} \label{P:uniform diameter 2}
There exists a $D > 0$ so that for any simplex $u \in \C^\Delta(S,z)$,
\[
\diam\big(\hull(\Gamma(u)) \cap N_1(\hull(G))\big) \leq D.
\]
\end{proposition}
\begin{proof}
Let $D$ be as in Corollary \ref{C:uniform diameter}.  If $u \in \widehat \C^\Delta(S,z)$, then the bound on diameter is precisely the conclusion of Corollary \ref{C:uniform diameter}.  If $u \in \C^\Delta(S,z) - \widehat \C^\Delta(S,z)$, then $u$ is preperipheral and so by our choice of horoball $\hull(\Gamma(u)) \cap N_1(\hull(G)) = \emptyset$.
\end{proof}

\bigskip
\noindent
Fix a vertex $u \in \C^0(S,z)$ and a point $x \in \hull(\Gamma(u)) \cap \hull(G)$.   A finite generating set for $G$
defines a word metric on $G$, but it is more convenient to use the metric
\[
\d_G(g,h) := \d_{\hull(G)}(g(x),h(x)) = \d_{\widetilde S}(g(x),h(x))
\]
which is quasiisometric to any such word metric by the Milnor-$\check{\mathrm{S}}$varc Lemma.
The following implies Theorem \ref{T:convex ccpt 1}.

\begin{theorem} \label{T:quasiisometric embedding}
The orbit map
$
G \to G \cdot u
$
given by $g \mapsto g \cdot u$ is a quasiisometric embedding into $\C(S,z)$.
\end{theorem}

\begin{proof}[Proof of Theorem \ref{T:convex ccpt 1} from Theorem \ref{T:quasiisometric embedding}.]
It was shown in \cite{kentleininger} (Theorem 1.3) and independently in \cite{hamenstadt} (Theorem 2.9) that a finitely
generated subgroup of the mapping class group is convex cocompact if and only if the orbit map to the curve complex is
a quasiisometric embedding.\end{proof}

\begin{proof}[Proof of Theorem \ref{T:quasiisometric embedding}.]
Write $\d_1$ for the metric on $\C^0(S,z)$ induced from the inclusion into $\C^1(S,z)$, see \S \ref{S:curve complex}.
We must find $K \geq 1$ and $C \geq 0$ so that for any $g \in G$, we have
\[
\frac{\d_G({\bf 1},g)}{K} - C \leq \d_1(u,g \cdot u) \leq K \d_G({\bf 1},g) + C.
\]

The upper bound follows from the fact that $\d_G$ is quasiisometric to the word metric, for which such an upper bound is an immediate consequence of the triangle inequality.  We assume that the constants $K$ and $C$ we choose for the lower bound also suffice for the upper bound.

We proceed to the proof of the lower bound.

Let $\tau \co  \widetilde S \to \hull(G)$ denote the closest point projection.  This is a contraction.  Moreover,
a well--known fact in hyperbolic geometry is that there exists an $R > 0$ so that if $\sigma$ is any geodesic segment outside $N_1(\hull(G))$ then
$\tau(\sigma)$ has length $\ell(\tau(\sigma)) \leq R$.

Next, suppose that $u'$ is a simplex in $\C(S,z)$ and $\sigma$ is a geodesic segment contained in $\hull(\Gamma(u'))$.
Since $\hull(\Gamma(u')) \cap N_1(\hull(G))$ is convex, it cuts $\sigma$ into at most three geodesic segments, at most one of which is contained in $\hull(\Gamma(u')) \cap N_1(\hull(G))$.  It follows that
\[
\ell(\tau(\sigma)) \leq 2R + D,
\]
where $D$ is as in Proposition \ref{P:uniform diameter 2}.

Now suppose that
\[
n = \d_1(u,g \cdot u)
\]
and connect $u$ to $g \cdot u$ by a geodesic edge path $[u_0,\ldots,u_n]$. It follows from the $\pi_1(S)$--equivariance of
$\hull \circ \Gamma$ (see Corollary \ref{C:complex to convex sets}) that $g(\hull(\Gamma(u))) = \hull(\Gamma(g \cdot u))$. We construct a
piecewise geodesic path
\[
\gamma \co [0,2n+1] \to \widetilde S
\]
connecting $x \in \hull(\Gamma(u))$ to $g(x) \in g(\hull(\Gamma(u))) = \hull(\Gamma(g \cdot u))$ as follows.

Consider the $2n+1$ simplices in the $1$--skeleton of $\C(S,z)$
\[
\begin{array}{rclrclrclrcl}
w_0 & = & \{u_0\}, & w_1 & = & \{u_0,u_1\},\\
w_2 & = & \{u_1\}, & w_3 & = & \{u_1,u_2\}, \\
w_4 & = & \{u_2\}, &     w_5 & = & \{u_2,u_3\}, \\
 & \vdots &    & & \vdots \\
w_{2n-2} & = & \{u_{n-1}\}, & w_{2n-1} & = &\{u_{n-1},u_n\},\\
w_{2n} & = & \{u_n\} \end{array}
\]
These have the property that
\[
w_{2j},w_{2j+2} \subset w_{2j+1}
\]
for every $j = 0,\ldots,n-1$.  Therefore, by Corollary \ref{C:complex to convex sets} it follows that
\[
\hull(\Gamma(w_{2j})), \hull(\Gamma(w_{2j+2})) \supseteq \hull(\Gamma(w_{2j+1}))
\]
for every $j = 0, \ldots,n-1$.
The key consequence is that
\[
\hull(\Gamma(w_k)) \cap \hull(\Gamma(w_{k+1})) \neq \emptyset
\]
for every $k = 0, \ldots,2n-1$.

From this it follows that we can define a path $\gamma \co [0,2n+1] \to \widetilde S$ with the following properties
\begin{itemize}
\item $\gamma(0) = x$
\item $\gamma(2n+1) = g(x)$
\item $\gamma([k,k+1]) \subset \hull(\Gamma(w_k))$ is a geodesic segment.
\end{itemize}

To see this, note that we have a chain of $2n+1$ convex sets, the first containing $x$ and the last containing $g(x)$.
Consecutive sets in the chain nontrivially intersect.  We therefore take $\gamma([0,1])$ to be the geodesic
segment from $x$ to \textit{any} point of the intersection $\hull(\Gamma(w_0)) \cap \hull(\Gamma(w_1))$.  Next take
$\gamma([1,2])$ to be the geodesic segment connecting this point to any point of $\hull(\Gamma(w_1)) \cap
\hull(\Gamma(w_2))$, and so on.  We continue in this way, ending with a geodesic segment $\gamma([2n,2n+1])$ connecting
the already determined point of $\hull(\Gamma(w_{2n-1})) \cap \hull(\Gamma(w_{2n}))$ to $g(x)$.  Convexity guarantees
the last property required.

Since the geodesic segment $\gamma([k,k+1])$ is contained in $\hull(\Gamma(w_k))$, we see that for every $k = 0, \ldots,2n$
\[
\ell(\tau(\gamma([k,k+1]))) \leq 2R +D
\]
Since $\gamma$ connects $x$ to $g(x)$, so does $\tau(\gamma)$, and its length bounds the distance from $x$ to $g(x)$.  Therefore we obtain
\[
\d_G({\bf 1},g) = \d_{\hull(G)}(x,g(x)) \leq \ell(\tau(\gamma)) \leq (2n+1)(2R+D)
\]
Isolating $n = \d_1(u,g \cdot u)$ in this inequality, we have
\[
\d_1(u,g \cdot u) = n \geq \frac{\d_G({\bf 1},g)}{2(2R+D)} - \frac{1}{2}
\]

Taking any $K \geq 2(2R+D)$ and $C \geq 1/2$ completes the proof.
\end{proof}

\section{Trees} \label{S:tree section}

Given a simplex $v \in \C^\Delta(S)$, there is an associated action of $\pi_1(S)$ on a tree $T_v$, namely, the Bass--Serre tree for the splitting of $\pi_1(S)$ determined by $v$.  We
refer the reader to \cite{Shalen} for a general introduction to actions on trees associated to codimension--$1$ submanifold.

In this section we prove
\begin{theorem} \label{T:main tree 1}
The fiber of $\Pi$ over a point $x$ in the interior of a simplex $v \subset \C(S)$ is $\pi_1(S)$--equivariantly
homeomorphic to the tree $T_v$ determined by $v$.
\end{theorem}

The fiber naturally inherits the structure of a metric simplicial tree from the point $x$, and as we vary this point in the base, we vary the metric trees continuously in the space of $\pi_1(S)$--trees, see \S \ref{S:metric}.  Harer defined a section of $\Pi\co \widehat \C(S,z) \to \C(S)$, and the metric on the trees can be used to parameterize a straight line deformation retraction to the image of the section.  This allows us to give an alternative proof of the following theorem.
\begin{corollary} [Harer]
If $\xi(S) \geq 1$, then with respect to the polyhedral topologies, $\C(S)$ is homotopy equivalent to $\widehat \C(S,z)$.
\end{corollary}
\noindent Hatcher and Vogtmann have given a simplified proof of this corollary, see \cite{hatchervogtmann}.

We begin with a discussion of the fibers of $\Pi$.

\subsection{Fibers} \label{S:fibers of Pi}

For any $x \in \C(S)$, the fiber $\mathcal F_x = \Pi^{-1}(x)$ can be naturally given the structure of a
simplicial complex $\mathcal F_x^\Delta$ so that each simplex is affinely embedded in a simplex of $\widehat \C(S,z)$.
Let $v$ be the unique simplex of $\C^\Delta(S)$ containing $x$ in its interior and let $\mathcal F_v^\Delta =
\Pi_\Delta^{-1}(v)$ endowed with the partial order obtained by restricting the partial order on $\widehat \C^\Delta(S,z)$
to $\mathcal F_v^\Delta$. We emphasize that $\Pi^{-1}_\Delta(v)$ is \textit{not} a simplicial subcomplex of $\widehat
\C^\Delta(S,z)$ (unless $v \in \C^0(S)$), but simply the partially ordered set of simplices sent by
$\Pi_\Delta$ to the simplex $v$.

Note that $\mathcal F_x^\Delta$ is $\pi_1(S)$--equivariantly order isomorphic to $\mathcal F_v^\Delta$: the isomorphism
is given by sending a simplex of $\mathcal F_x^\Delta$ to the smallest simplex of $\mathcal F_v^\Delta$ containing it.

\subsection{The trees} \label{S:trees and hulls}

We now recall one construction of the tree $T_v$ by constructing its underlying simplicial complex $T_v^\Delta$.
The simplicial tree $T_v^\Delta$ is the tree dual to the preimage $p^{-1}([v]) \subset
\widetilde S$.  More precisely, the vertices of $T_v^\Delta$ are in a one-to-one correspondence with components of
$\widetilde S - p^{-1}([v])$ with two vertices joined by an edge if the closures of the corresponding components
nontrivially intersect.

The edge and vertex stabilizers of $T^\Delta_v$ in $\pi_1(S)$ are precisely the stabilizers of the components of $p^{-1}([v])$ and of $\widetilde S \setminus p^{-1}([v])$, respectively.  By construction, the quotients of these by their stabilizers inject as components of $[v]$ and $S \setminus [v]$.  It follows from the discussion in Section \ref{S:group to hull} that these subgroups are all contained in $\D$ (falling into the two first cases (i) and (ii)).  Indeed, setting $\D_v$ to be the set of all edge and vertex stabilizers of $T^\Delta_v$ we have
\[
\D_v = \{ \Gamma(u) \, | \, \Pi(u) = v \} \subset \D.
\]

The group $\pi_1(S)$ acts on $\D_v$ by conjugation (the restriction of the action on $\D$).  Notice that
the stabilizers of two distinct vertices are distinct, and similarly for the edge stabilizers.  Thus, the
stabilizer determines the simplex of the tree. Therefore, since (in this setting) the stabilizer of a vertex properly contains the
stabilizer of any edge having it as an endpoint, we see that the simplicial complex $T_v^\Delta$ is $\pi_1(S)$--equivariantly
 reverse--order  isomorphic to $\D_v$. We record this as a proposition.

\begin{proposition} \label{P:tree groups}
For any simplex $v \in \C^\Delta(S)$, $T_v^\Delta$ is $\pi_1(S)$--equivariantly  reverse--order  isomorphic to $\D_v$.\qed
\end{proposition}

The map $\Gamma$ restricts to a map
\[
\Gamma_v \co   \mathcal F_v^\Delta \to \D_v.
\]
given by $\Gamma_v = \Gamma|_{\mathcal F_v^\Delta}$.
\begin{proposition} \label{P:fiber groups}
For any simplex $v \in \C^\Delta(S)$, $\Gamma_v$ is a reverse--order bijection.
\end{proposition}

The proof will require the following lemma.
\begin{lemma} \label{L:geomtop}
Suppose $(Y,z),(Y',z) \subset (S,z)$ are compact $\pi_1$--injective subsurfaces with $Y,Y' \ncong D^2$ and $\pi_1(Y,z) = \pi_1(Y',z) < \pi_1(S,z)$ proper.
Then there is an isotopy fixing $z$ taking $Y$ to $Y'$.
\end{lemma}
\noindent
It is easy to see that there is an isotopy taking $Y$ to $Y'$ if we do not require such an isotopy to fix $z$.
The proof of the lemma is an exercise in geometric topology, which we sketch for completeness.
\begin{proof}[Sketch]
With our chosen isomorphism $\pi_1(S) \cong \pi_1(S,z)$, we view $\pi_1(Y,z) = \pi_1(Y',z)$ as a subgroup $G < \pi_1(S)$ of the covering group of $p \co \widetilde S \to S$.
First, by an isotopy fixing $z$ we may assume that $Y$ and $Y'$ meet transversely (equivalently, their boundaries meet transversely).
By further isotopy fixing $z$ we may assume that $\partial Y$ meets $\partial Y'$ in the fewest possible number of points.
If these boundaries are disjoint, then we easily find annuli disjoint from $z$ which we can use to produce an isotopy taking $Y$ to $Y'$.

Now suppose that the boundaries intersect.  We describe how to find a bigon not containing $z$ bounded by arcs of $\partial Y$ and $\partial Y'$ which can be used to reduce the number of intersection points, producing a contradiction.  To this end, consider $\widetilde Y$ and $\widetilde Y'$, the $G$--invariant components of $p^{-1}(Y)$ and $p^{-1}(Y')$, respectively.  By assumption, $\widetilde Y$ and $\widetilde Y'$ contain exactly the same subset of $p^{-1}(z)$, namely all $G$--translates of the chosen basepoint $\widetilde z \in p^{-1}(z)$ defining the isomorphism $\pi_1(S) \cong \pi_1(S,z)$.  The same is true for the translate of $\widetilde Y$ and $\widetilde Y'$ by any element of $\pi_1(S)$.  Said differently, $p^{-1}(\partial Y)$ and $p^{-1}(\partial Y')$ define exactly the same partition of $p^{-1}(z)$.  An  innermost  bigon in $\widetilde S$ bounded by arcs of $p^{-1}(\partial Y)$ and $p^{-1}(\partial Y')$ projects to the desired bigon in $S$.
\end{proof}

\begin{proof}[Proof of Proposition \ref{P:fiber groups}]
As we have already noted, $\Gamma_v$ is an order--reversing surjection.
We must show that $\Gamma_v$ is injective.

To see this, we suppose $\Gamma(u) = \Gamma(u')$ with $u,u' \in \mathcal F_v^\Delta$.  We need to show that $u = u'$.
From the definition of $\Gamma$ we have $\pi_1(Y(u),z) = \pi_1(Y(u'),z)$ (recall that we must realize $u$ and $u'$ by multicurves to make sense of $Y(u)$ and $Y(u')$).
By Lemma \ref{L:geomtop} there is an isotopy fixing $z$ taking $Y(u)$ to $Y(u')$.  This proves that $\partial_0 Y(u) = \partial_0 Y(u')$.

Now let $u_0 \in u$ and $u_0' \in u'$ be vertices \textit{not} in $u$ and $u'$, respectively, for which $\Pi(u_0) = \Pi(u_0')$ (if there are no such vertices, then $u = \partial_0 Y(u) = \partial_0 Y(u') = u'$ and we are done).  We must show that $u_0 = u_0'$.
This is another argument in geometric topology.  By an isotopy fixing $Y(u) = Y(u')$ first assume $u_0$ and $u_0'$ intersect transversely in the fewest possible number of points (fewest among all curves isotopic by an isotopy fixing $Y(u) = Y(u')$).  If $u_0$ and $u_0'$ are disjoint, then because they become isotopic after forgetting $z$, they must cobound an annulus in $S$.  This annulus cannot possibly contain $Y(u) = Y(u')$, so it may be used to produce an isotopy fixing $Y(u) = Y(u')$ taking $u_0$ to $u_0'$ as required.
If they are not disjoint, then there is a bigon in $S$ bounded by arcs of these curves.  Of course this bigon cannot contain $Y(u) = Y(u')$, and so it can be used to produce an isotopy fixing $Y(u) = Y(u')$ and reducing the number of intersection points of $u_0$ and $u_0'$, as required.
Therefore $u_0 = u_0'$, so $u = u'$, and $\Gamma_v$ is injective.
\end{proof}

\begin{proof}[Proof of Theorem \ref{T:main tree 1}.]
Since $\Gamma_v$ is a reverse--order bijection and since $T_v^\Delta$ is reverse--order isomorphic to $\D_v$, it follows that we have an order--preserving bijection from $T^\Delta_v$ to $\mathcal F_x^\Delta$.  Since these are both abstract $1$--dimensional simplicial complexes, they must be isomorphic.
\end{proof}

\subsubsection{Variation of metrics}\label{S:metric}

The simplex $v$ is naturally a space of transverse measures on its underlying $1$--manifold: namely, a point in $v$ is a convex combination of the unit weights on the components of this $1$--manifold.

The tree $\mathcal F_x$ thus inherits a metric varying continuously in $x$ by assigning to each edge the weight of the dual curve in $x$ to which it corresponds.
This metric is the same as the path metric induced by the inclusion of $\mathcal F_x$ into $\widehat \C(S,z)$.

\subsection{Harer's section and a deformation retraction}

There are many ways to construct a section of the map $\Pi \co \widehat \C(S,z) \to \C(S)$.
Harer describes one as follows.  First, the union of the geodesic representatives
\[
\bigcup_{v \in \C^0(S)} [v]
\]
has measure zero.  We chose our basepoint $z$ to be any point outside this union.
A section is then given by $v \mapsto [v]$.  This makes sense because $[v]$ lies in $S \setminus \{z \}$ and geodesics minimize intersection between pairs of curves.

Let us denote the image of this section by $\C'(S)$.
The map $\Pi \co \widehat \C(S,z) \to \C(S)$ composed with this section gives a map $\Pi' \co \widehat \C(S,z) \to \C'(S)$.
The fibers of $\Pi'$ are precisely the fibers of $\Pi$ and are therefore metric trees.  
The section provides a preferred basepoint in each, the intersection with $\C'(S)$.

Any metric tree $T$ with a preferred basepoint $x$ admits a ``straight line'' deformation retraction to the basepoint
\[
H \co  T \times [0,1] \to T
\]
defined by setting $H(y,t)$ to be the unique point of the geodesic segment $[x,y]$ for which $d(x,H(y,t)) = (1-t)d(x,y)$.

This determines a map
\[
H \co  \widehat C(S,z) \times [0,1] \to \widehat \C(S,z)
\]
defined on each of the pointed trees by the procedure just described.

This map is
not continuous with respect to the metric topology on $\widehat \C(S,z)$.
The idea  is that for a
vertex $u$ in $\widehat \C(S,z)$ and any $0 < \epsilon < 1$ one can find a point $x$ within $\epsilon$ of $u$ for which
$H(\{x\} \times [0,1-\epsilon])$ is always within $\epsilon$ of $u$, while $H(u,t)$ is making progress toward $\C'(S)$.
In particular, one can construct sequences $\{x_n\}$ with $x_n \to u$ and $H(x_n,1/2) \to u$ but $H(u,1/2)$ far from
$u$.

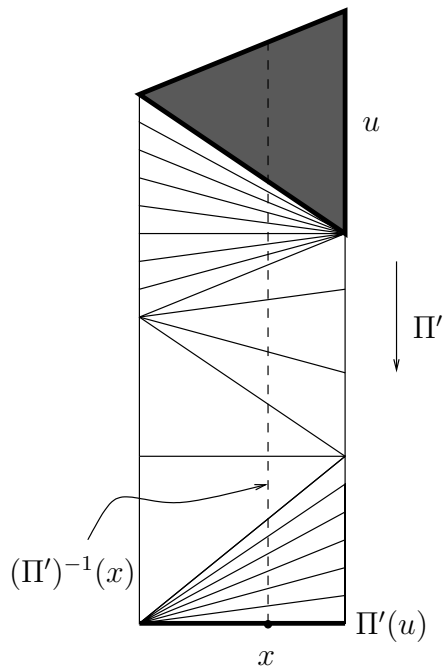
\begin{figure}[h!]
\begin{center}
\input{retract.pstex_t}
\end{center}
\caption{The deformation retraction above a simplex $\Pi'(u)$.}
\label{F:retract}
\end{figure}

However, the polyhedral topology is more natural from the perspective of algebraic topology, and here the map $H$ is continuous.
\begin{proposition}
The map $H$ is continuous with respect to the polyhedral topologies, and hence is a deformation retraction.
\end{proposition}
\begin{proof}
With respect to the polyhedral topologies, it suffices to show that the restriction of $H$ to $u \times [0,1]$ is continuous for each $u \in \widehat \C^\Delta(S,z)$.  The set
\[
H(u \times [0,1])
\]
is a finite subcomplex $X \subset \widehat \C(S,z)$: it is the union of paths in the fiber of $\Pi'$ from $u$ to $\Pi'(u)$, and there are only finitely many combinatorial types of these paths, one for each face of $\Pi'(u)$.  Thus, the restriction becomes a map
\[
H|_{u \times [0,1]}  \co  u \times [0,1] \to X
\]
Since $X$ is a finite simplicial complex, it is easy to check that $H|_{u \times [0,1]}$ is continuous. Figure
\ref{F:retract} gives a cartoon of the general situation where the euclidean simplices have been distorted affinely.
\end{proof}

\bibliographystyle{plain}
\bibliography{tree}

\bigskip

\noindent Department of Mathematics, Brown University, Providence, RI 02912 \newline \noindent
\texttt{rkent@math.brown.edu}

\bigskip

\noindent Department of Mathematics, University of Illinois, Urbana--Champaign, IL 61801 \newline \noindent  \texttt{clein@math.uiuc.edu}

\bigskip

\noindent Department of Mathematics, University of Warwick, Coventry CV4 7AL UK
\newline \noindent \texttt{s.schleimer@warwick.ac.uk }

\end{document}

%% file: subsurface.pstex_t
\begin{picture}(0,0)%
\includegraphics{subsurface.pstex}%
\end{picture}%
\setlength{\unitlength}{3947sp}%
\begingroup\makeatletter\ifx\SetFigFont\undefined%
\gdef\SetFigFont#1#2#3#4#5{%
  \reset@font\fontsize{#1}{#2pt}%
  \fontfamily{#3}\fontseries{#4}\fontshape{#5}%
  \selectfont}%
\fi\endgroup%
\begin{picture}(5601,3728)(1284,-5753)
\put(2664,-2799){\makebox(0,0)[lb]{\smash{{\SetFigFont{11}{13.2}{\familydefault}{\mddefault}{\updefault}{\color[rgb]{0,0,0}$z$}%
}}}}
\put(5633,-2851){\makebox(0,0)[lb]{\smash{{\SetFigFont{11}{13.2}{\familydefault}{\mddefault}{\updefault}{\color[rgb]{0,0,0}$z$}%
}}}}
\put(4393,-4562){\makebox(0,0)[lb]{\smash{{\SetFigFont{11}{13.2}{\familydefault}{\mddefault}{\updefault}{\color[rgb]{0,0,0}$z$}%
}}}}
\put(2174,-3878){\makebox(0,0)[lb]{\smash{{\SetFigFont{11}{13.2}{\familydefault}{\mddefault}{\updefault}{\color[rgb]{0,0,0}$u$}%
}}}}
\put(6067,-3483){\makebox(0,0)[lb]{\smash{{\SetFigFont{11}{13.2}{\familydefault}{\mddefault}{\updefault}{\color[rgb]{0,0,0}$Y(u)$}%
}}}}
\put(3887,-5699){\makebox(0,0)[lb]{\smash{{\SetFigFont{11}{13.2}{\familydefault}{\mddefault}{\updefault}{\color[rgb]{0,0,0}$\partial_0 Y(u)$}%
}}}}
\end{picture}%

%% file: retract.pstex_t
\begin{picture}(0,0)%
\includegraphics{retract.pstex}%
\end{picture}%
\setlength{\unitlength}{3947sp}%
\begingroup\makeatletter\ifx\SetFigFont\undefined%
\gdef\SetFigFont#1#2#3#4#5{%
  \reset@font\fontsize{#1}{#2pt}%
  \fontfamily{#3}\fontseries{#4}\fontshape{#5}%
  \selectfont}%
\fi\endgroup%
\begin{picture}(2535,4197)(2506,-5300)
\put(4726,-1891){\makebox(0,0)[lb]{\smash{{\SetFigFont{12}{14.4}{\familydefault}{\mddefault}{\updefault}{\color[rgb]{0,0,0}$u$}%
}}}}
\put(4681,-5056){\makebox(0,0)[lb]{\smash{{\SetFigFont{12}{14.4}{\familydefault}{\mddefault}{\updefault}{\color[rgb]{0,0,0}$\Pi'(u)$}%
}}}}
\put(5041,-3196){\makebox(0,0)[lb]{\smash{{\SetFigFont{12}{14.4}{\familydefault}{\mddefault}{\updefault}{\color[rgb]{0,0,0}$\Pi'$}%
}}}}
\put(4066,-5251){\makebox(0,0)[lb]{\smash{{\SetFigFont{12}{14.4}{\familydefault}{\mddefault}{\updefault}{\color[rgb]{0,0,0}$x$}%
}}}}
\put(2506,-4726){\makebox(0,0)[lb]{\smash{{\SetFigFont{12}{14.4}{\familydefault}{\mddefault}{\updefault}{\color[rgb]{0,0,0}$(\Pi')^{-1}(x)$}%
}}}}
\end{picture}%